\documentclass[11pt,reqno]{amsart}
\usepackage{amsmath}
\usepackage{cases}
\usepackage{mathrsfs}
\usepackage{bbm}
\usepackage{amssymb}
\usepackage{amscd}
\usepackage{amsfonts,latexsym,amsmath,amsthm,amsxtra,mathdots,amssymb,latexsym,mathabx}
\usepackage[all,cmtip]{xy}
\RequirePackage{amsmath} \RequirePackage{amssymb}
\usepackage{color}
\usepackage{colordvi}
\usepackage{multicol}
\usepackage{hyperref}
\usepackage{mathtools}

\usepackage[margin=1in]{geometry}
\usepackage{xcolor}
\hypersetup{
	colorlinks,
	linkcolor={red!50!black},
	citecolor={blue!50!black},
	urlcolor={blue!80!black}
}
\usepackage{cite}

\newcommand{\bea}{\begin{eqnarray}}
	\newcommand{\eea}{\end{eqnarray}}
\newcommand{\bna}{\begin{eqnarray*}}
	\newcommand{\ena}{\end{eqnarray*}}

\numberwithin{equation}{section}

\theoremstyle{plain}
\newtheorem*{Theorem A}{Theorem A}
\newtheorem*{Theorem B}{Theorem B}

\newtheorem{lemma}{Lemma}[section]
\newtheorem{theorem}{Theorem}[section]
\newtheorem{corollary}{Corollary}[section]
\newtheorem{proposition}[lemma]{Proposition}

\theoremstyle{definition}

\newtheorem{remark}{Remark}

\begin{document}
	
	\title
	[Note on shifted primes with large prime factors]
	{Note on shifted primes with large prime factors}
	
	\author
	[Y. Ding \quad $\&$ \quad   Z. Wang] {Yuchen Ding \quad \text{and} \quad Zhiwei Wang}
	
	\address{(Yuchen Ding) $^1$School of Mathematics,  Yangzhou University, Yangzhou 225002, Jiangsu, China; 
		$^2$Alfr\'ed R\'enyi Institute of Mathematics (LERN), Re\'altanoda utca 13-15, H-1053 Budapest, Hungary}
	\email{\tt ycding@yzu.edu.cn}
	
	\address{(Zhiwei Wang)  $^1$School of Mathematics,  Shandong University, Jinan 250100, Shandong, China; $^2$State Key Laboratory of Cryptography and Digital Economy Security, Shandong 
		University, Jinan 250100, Shandong, China
	}
	\email{\tt zhiwei.wang@sdu.edu.cn}

	\subjclass[2010]{Primary 11N05; Secondary 11N36}
	
	\keywords{shifted primes, largest prime factor, linear sieve, primes in arithmetic progressions}

	\begin{abstract}
		We denote by $P^+(n)$ the largest prime factor of the integer $n$. 
		In 1935, Erd\H os studied the quantity $T_c(x)$ defined by 
		$$
		T_c(x)=\big|\big\{p\le x:  P^+(p-1)\ge p^c\big\}\big|,
		$$
		and he proved 
		$$
		\limsup_{x\rightarrow \infty}\frac{T_c(x)}{\pi(x)}\rightarrow 0, \quad \text{as~}c\rightarrow 1.
		$$
		Recently, Ding gave a  quantitative form of Erd\H os' result, showing that  
		$$
		\limsup_{x\rightarrow \infty}\frac{T_c(x)}{\pi(x)}\le 8\big(c^{-1}-1\big).
		$$
		holds for $8/9< c<1$. In this paper, we improve Ding's upper bound to
		$$
		\limsup_{x\rightarrow \infty}\frac{T_c(x)}{\pi(x)}\le -\frac{7}{2}\log c
		$$
		for $e^{-\frac{2}{7}}< c<1$.
	\end{abstract}
	
	\maketitle
	
	\section{Introduction}
	We denote by $P^+(n)$ the largest prime factor of an integer $n$,
	with the convention that $P^+(1)=1$. The study of largest prime factor of shifted prime $P^+(p+a), a\in \mathbb{Z}$ is of significant importance. First,  the infinitude of primes $p$ with $P^+(p+2)>p$ is equivalent to the twin prime conjecture, which is one of the most well-known open problems in number theory; second, an unexpected connection  between large value of $P^+(p-1)$ and the first case of Fermat's last theorem, was established by Adleman and Heath-Brown \cite{AHB}, and Fouvry \cite{Fouvry}. Last but not least, small values of $P^+(p+a)$ plays an important role in cryptography, such as Pollard's $p-1$ algorithm and Williams' $p+1$ algorithm.
	
	In this article, we study the quantity $T_c(x)$ defined by
	$$
	T_c(x):=\big|\big\{p\le x:  P^+(p-1)\ge p^c\big\}\big|,
	$$
	where $0<c<1$.
	
	For small values of $c\le1/2$, here we just list several results on $T_c(x)$. 
	In 2015, Luca, Menares, and Pizarro-Madariaga \cite{LMP} showed that
	\begin{align}\label{eq-1-2}
		\liminf_{x\rightarrow\infty}\frac{T_c(x)}{\pi(x)}\ge 1-c
	\end{align}
 holds for $1/4\le c\le 1/2$. 
	Later, Chen and Chen \cite{CC} extended the range of $c$ to $(0,1/2)$ for (\ref{eq-1-2}).
	In 2018, Feng and Wu \cite{FW} improved the lower bound in (\ref{eq-1-2}) for $0<c< 0.3517\ldots$ by showing
	$$
	\liminf_{x\rightarrow\infty}\frac{T_c(x)}{\pi(x)}\ge 1-4\int_{1/c-1}^{1/c}\frac{\rho(t)}{t}dt,
	$$
	where $\rho(u)$ is the Dickman function. 
	 Later, the above lower bound was further improved by Liu, Wu and Xi \cite{LWX} to
	$$
	\liminf_{x\rightarrow\infty}\frac{T_c(x)}{\pi(x)}\ge 1-4\rho(1/c)
	$$
	provided $0< c<0.3734\ldots$. Finally, we refer the readers to e.g. Banks and Shparlinski \cite{Banks}, Wu \cite{Wu-1,Wu-2} for more results related to this topic.

	Now we turn our attention to the case when $1/2<c\le 1$ is large. It seems this case is more difficult and more attractive, considering its close connection with the twin prime conjecture and Fermat's last theorem as we have mentioned above.
	On the one hand, one may want to seek values of $c$ as large as possible such that
	\begin{align}\label{eq-1-1}
		\liminf_{x\rightarrow\infty}\frac{T_c(x)}{\pi(x)}>0.
	\end{align}	
	In 1969, Goldfeld \cite{Goldfeld} proved that \eqref{eq-1-1} holds for some $c>1/2$ by employing the Bombieri-Vinogradov theorem as well as the Brun-Titchmarsh inequality, and he also stated that his arguments is also applicable for  $c\approx 7/12$. Since then,
there has been a number of improvements on the lower bound of values of $c$ in (\ref{eq-1-1}), see, e.g. Motohashi \cite{Motohashi}, Hooley \cite{Hoo1,Hoo2}, Deshouillers and Iwaniec \cite{Iwaniec}, Fouvry \cite{Fouvry}, Baker and Harman \cite{BH}. The best record of $c$ up to now is $0.677$, obtained by Baker and Harman \cite{BH-2}. 
	
	On the other hand, one may want to have a non-trivial upper bound for  $T_c(x)$ for  some fixed large $c$, which is the main topic we shall study in this paper. In fact, we even conditionally know the asymptotic formula for  $T_c(x)$: 
	\begin{align}\label{asymp}
		\lim_{x\rightarrow\infty}\frac{T_c(x)}{\pi(x)} =1-\rho\left(\frac{1}{c}\right)
	\end{align}
	holds for $0<c<1$ under the assumption of the Elliott--Halberstam conjecture, see e.g. the works of Pomerance \cite{Pomerance}, Granville \cite{Granville}, Wang \cite{Wang} and Wu \cite{Wu}. Motivated by this, in 2023, Ding \cite{Ding} proved unconditionally that
	\begin{align}\label{eq-1-4}
		\limsup_{x\rightarrow\infty}\frac{T_c(x)}{\pi(x)}<1/2
	\end{align}
	for some absolute constant $c<1$. 
	The proof of (\ref{eq-1-4}) by Ding is based on the following result of Brun-Titchmarsh type due to Wu \cite[Lemma 2.2]{Wu}. 
		\begin{proposition}\label{pro1}
		For $(a,m)=1$, let $\pi(x;m,a)$ denotes the number of primes $p\le x$ such that $p\equiv a\pmod{m}$.There exist two functions $K_2(\theta)>K_1(\theta)>0$, defined on the interval $(0,17/32)$ such that for each fixed  $A>0$, and sufficiently large $Q=x^\theta$, the inequalities
		$$K_1(\theta)\frac{\pi(x)}{\varphi(m)}\le \pi(x;m,1)\le K_2(\theta)\frac{\pi(x)}{\varphi(m)}$$
		hold for all integers $m\in(Q,2Q]$ with at most $O\left(Q(\log Q)^{-A}\right)$ exceptions, where the implied constant depends only on $A$ and $\theta$. Moreover, for any fixed $\varepsilon>0$, these functions can be chosen to satisfy the following properties:\\
		$\bullet$ $K_1(\theta)$ is monotonic decreasing, and $K_2(\theta)$ is monotonic increasing.\\
		$\bullet$ $K_1(1/2)=1-\varepsilon$ and $K_2(1/2)=1+\varepsilon$.
	\end{proposition}
	
	The constant $c$ in (\ref{eq-1-4}) could further be specified by explicit values of $K_1(\theta)$ in Proposition \ref{pro1}, for example, one has $K_1(\theta)\ge 0.16$ for $1/2\le\theta\le 13/25$ \cite[Theorem 1]{BH}. Using the method of Ding \cite{Ding} as well as the explicit values of $K_1(\theta)$, Xinyue Zang (private communication) obtained that
\begin{align}\label{eq-1-5}
	\limsup_{x\rightarrow\infty}\frac{T_c(x)}{\pi(x)}\le \frac{0.496875}{c}<\frac{1}{2}, \quad \text{for~} 0.993375<c<1.
\end{align}

However, there are earlier results related to (\ref{eq-1-4}) and (\ref{eq-1-5}).  Actually, as indicated by the proof of a former result of Erd\H os \cite[from line -6, page 212 to line 4, page 213]{erdos}, as early as 1935, one could already conclude that 
\begin{align}\label{eq-1-6}
	\limsup_{x\rightarrow\infty}\frac{T_c(x)}{\pi(x)}\rightarrow 0, \quad \text{as}~c\rightarrow1,
\end{align}
which just corresponds to the asymptotic \eqref{asymp} since $1-\rho(1/c)\to 0$ as $c\to 1$. Clearly, (\ref{eq-1-4}) is now a simple corollary of (\ref{eq-1-6}). Very recently, Ding \cite{Ding-2} obtained a quantitative form of Erd\H os' result  (\ref{eq-1-6}):
\begin{align}\label{eq-1-7}
	\limsup_{x\rightarrow\infty}\frac{T_c(x)}{\pi(x)}\le 8\big(c^{-1}-1\big)
\end{align}
for $8/9< c<1$. By (\ref{eq-1-7}) one notes easily that 
\begin{align}\label{eq-add-1}
	\limsup_{x\rightarrow\infty}\frac{T_c(x)}{\pi(x)}<\frac{1}{2}
\end{align}
for any $16/17<c<1$ which improved the numerical values of (\ref{eq-1-5}). It should be mentioned that almost the same time as Ding's result \eqref{eq-1-7},  Bharadwaj and Rodgers \cite{BhaRod} independently obtained the same
result \eqref{eq-1-6} with a general form in probabilistic language.\footnote{All of the authors (Ding, Bharadwaj and Rodgers) were unaware of Erd\H os' result at an earlier time.}  Erd\H os' result (\ref{eq-1-6}) is an application of Brun's method, while the proof of (\ref{eq-1-7}) is mainly based on the following quantitative version of Selberg's upper bound sieve. 

	\begin{proposition}\label{pro2}\cite[Theorem 5.7]{Halberstam}
		Let $g$ be a natural number, and let $a_i,b_i~(i=1,2,\cdot\cdot\cdot,g)$ be integers satisfying 
		$$
		E:=\prod_{i=1}^{g}a_i\prod_{1\le r<s\le g}(a_rb_s-a_sb_r)\neq0.
		$$
		Let $\varrho(p)$ denote the number of solutions $n\pmod{p}$ to the congruence
		$$
		\prod_{i=1}^{g}(a_in+b_i)\equiv 0\pmod{p},
		$$
		and suppose that
		$$
		\varrho(p)<p \quad \text{for~all~} p.
		$$
		Let $y$ and $z$ be real numbers satisfying
		$
		1<y\le z.
		$ 
		Then we have
		\begin{multline*}
			\big|\big\{n:z-y<n\le z, a_in+b_i~\text{prime~for~}i=1,2,\cdot\cdot\cdot,g\big\}\big|\\
			\le 2^gg!\prod_{p}\left(1-\frac{\varrho(p)-1}{p-1}\right)\left(1-\frac{1}{p}\right)^{-g+1}\!\!\!\frac{y}{\log ^gy} \left(1+O\left(\frac{\log\log3y+\log\log 3|E|}{\log y}\right)\right),
		\end{multline*}
		where the constant implied by the $O$-symbol depends at most on $g$.
	\end{proposition}
	
	For the proof of (\ref{eq-1-7}), one used Proposition \ref{pro2} in the particular case $g=2$. Hence, the constant factor $8$ in (\ref{eq-1-7}) comes from the identity
	$
	2^gg!=8.
	$
	
	In this article, we shall give a further improvement of (\ref{eq-1-7}) with two new ingredients: the first is the employment of linear sieve to the prime variable sequence, instead of integer variable polynomial combining with the two dimensional sieve; the second one is that when dealing with the error term coming from linear sieve, we apply a theorem of Bombieri-Friedlander-Iwaniec type with level of distribution $x^{4/7-\varepsilon}$ instead of the classical level $x^{1/2-\varepsilon}$.
	
	Our main result is stated as follows.
	
	\begin{theorem}\label{thm1}
		For any $e^{-\frac{2}{7}}< c<1$ we have
		\begin{align*}
			\limsup_{x\rightarrow\infty}\frac{T_c(x)}{\pi(x)}\le -\frac{7}{2}\log c.
		\end{align*}
	\end{theorem}
	
	\begin{remark}\label{remark1}
		Theorem \ref{thm1} provides a nontrivial upper bound of $T_c(x)$ for any $e^{-\frac{2}{7}}< c<1$. Here the lower bound of $c$ is approximately
		$
		e^{-\frac{2}{7}}=0.75147\ldots,
		$ 
		which could be compared to $8/9=0.88888\ldots$ in (\ref{eq-1-7}). Thus, Theorem \ref{thm1} extends the range of $c$ in (\ref{eq-1-7}). Furthermore, one may see easily that
		$$
		-\frac{7}{2}\log c<-\frac{7}{2}\big(c^{-1}-1\big)<8\big(c^{-1}-1\big)
		$$
		for any $8/9< c<1$ and hence Theorem \ref{thm1} also improves the upper bound of $T_c(x)$ in (\ref{eq-1-7}).
	\end{remark}
	
\begin{remark}\label{remark-new-1}
Recall that $\rho(u)=1$ for $0\le u\le 1$ and 
$
\rho(u)=1-\log u
$
for $1\le u\le 2$ from the recursion formula (see e.g. \cite[Eq. (7.6)]{Montgomery}) on Dickman's function. Thus, by (\ref{asymp}) we have conditionally for $1/2<c<1$
$$
\lim_{x\rightarrow\infty}\frac{T_c(x)}{\pi(x)}=1-\rho\left(\frac{1}{c}\right)=-\log c
$$
under the assumption of Elliott--Halberstam conjecture, which seems  far out of reach at present. Our unconditionally bound $-\frac{7}{2}\log c$ in Theorem \ref{thm1} is 3.5 times that of the expected asymptotic density.
			\end{remark}	
	
	The following corollary of Theorem \ref{thm1} improves (\ref{eq-1-4}), (\ref{eq-1-5}) and (\ref{eq-add-1}) considerably.
	\begin{corollary}\label{corollary}
		For any $c>e^{-\frac{1}{7}}$ we have $\limsup_{x\rightarrow\infty}T_c(x)/\pi(x)<1/2$.
	\end{corollary}
	\begin{remark}\label{remark2}
		The numerical value of $e^{-\frac{1}{7}}$ is $\approx 0.86687$. 
		According to the conditional asymptotic formula \eqref{asymp}, one should expect $\limsup_{x\rightarrow\infty}T_c(x)/\pi(x)<1/2$ for any $c>e^{-\frac{1}{2}}\approx 0.60653$.
Corollary \ref{corollary} makes some further progress toward this direction.
	\end{remark}
	
	\section{Fundamental lemmas}
	
	 Let $\mu(n)$ be the M\"obius function. Let $\mathcal{A}$ be a finite sequence of positive integers and let $\mathcal{P}$ be a subset of primes. For any $z\ge 2$, we set 
	$$P(z)=\prod_{\substack{ p\le z\\ p\in\mathcal{P}}}p.$$ 
	The first lemma is due to Iwaniec, see \cite{Iwaniec1} or \cite{Iwaniec2}.
	
	\begin{lemma}\label{lem:sieve}
		Let $D\geqslant z\geqslant 2 $ and $L>1$. Put $\lambda_1^{\pm}=1$ and for all squarefree numbers $d=p_1\cdots p_r,\, p_1>\cdots >p_r,\, r\geq 1$, define
		the upper and lower bound sieve weights of level $D$ respectively
		\begin{align*}
	\	\qquad	\lambda_d^{+}=\left\{ 
			\begin{array}{ll}
				(-1)^r &  \ \  {\rm if}\ p_{_1}\cdots p_{_{2l}}p_{_{2l-1}}^3<D \ {\rm for\ all\ } 0\leq l\leq (r-1)/2,
				\\\noalign{\vskip 2mm}
				0 &  \ \ {\rm otherwise},
			\end{array}
			\right.
		\end{align*}
and		
		\begin{align*}
			\lambda_d^{-} =\left\{
			\begin{array}{ll}
				(-1)^r &  \ \  {\rm if}\ p_{_1}\cdots p_{_{2l-1}}p_{_{2l}}^3<D \ {\rm for\ all\ } 0\leq l\leq r/2,
				\\\noalign{\vskip 2mm}
				0 &  \ \ {\rm otherwise}.
			\end{array}
			\right.
		\end{align*}
Then we have
		$$\lambda_d^- * \mathbf{1}\leqslant \mu * \mathbf{1} \leqslant \lambda_d^+ * \mathbf{1}$$ and 
		\begin{align*}
			\sum_{d \mid P(z)} \lambda^+_{d}\frac{ \omega(d)}{d} &\leqslant \prod_{p \leqslant z \atop p \in \mathcal{P}}\left(1-\frac{\omega(p)}{p}\right)
			\left\{F(s)+O\left(\frac{e^{\sqrt{L}-s}}{(\log D)^{1 / 3}}\right)\right\},
			\\\noalign{\vskip 1.6mm}
			\sum_{d \mid P(z)} \lambda^-_{d}\frac{ \omega(d)}{d} &\geqslant \prod_{p \leqslant z \atop p \in \mathcal{P}}\left(1-\frac{\omega(p)}{p}\right)
			\left\{f(s)-O\left(\frac{e^{\sqrt{L}-s}}{(\log D)^{1 / 3}}\right)\right\}
		\end{align*}
		uniformly for all multiplicative functions $w$ satisfying
		\begin{align*}
			& {\rm{(i)}} \quad 0<\omega(p)<p \quad(p \in \mathcal{P}),
			\\\noalign{\vskip 1.6mm}
			& {\rm{(ii)}} \prod_{\substack{u<p \leqslant v\\[1mm] p \in \mathcal{P}}}\left(1-\frac{\omega(p)}{p}\right)^{-1} \leqslant \frac{\log v}{\log u}\left(1+\frac{L}{\log u}\right) \quad(2 \leqslant u \leqslant v \leqslant z).
		\end{align*}
		where $s=\log D/\log z$ and $F(s), f(s)$  are defined by the continuous solutions to the system
		\begin{align*}
			\left\{
			\begin{array}{ll}
				sF(s)=2\mathrm{e}^{\gamma} &  \ \  1\leq s\leq 2,
				\\\noalign{\vskip 2mm}
				sf(s)=0 &  \ \ 0<s\leq 2,
				\\\noalign{\vskip 2mm}
				(sF(s))'=f(s-1) & \ \ s>2,
				\\\noalign{\vskip 2mm}
				(sf(s))'=F(s-1) & \ \ s>2.
			\end{array}
			\right.
		\end{align*}
		Here $\gamma$ is the Euler constant.
		\end{lemma}

	\vskip 3mm
The second lemma is a convolution generalization  theorem of Bombieri-Friedlander-Iwaniec \cite{BFI86} with well-factorable weights. An arithmetic function $\lambda(q)$ is called well factorable of level $Q$ if for any $Q=Q_1Q_2$, $Q_1, Q_2\ge 1$, there exist two functions $\lambda_1$ and $\lambda_2$ supported in $[1, Q_1]$ and $[1, Q_2]$ respectively such that
	$$
	|\lambda_1|\le 1, \quad |\lambda_2|\le 1 \quad \text{and} \quad \lambda=\lambda_1\ast \lambda_2.
	$$
Next for given integers $\ell, a, q $ with $(a,q)=1$, we define the convolution sum 
	$$
	\pi(y;\ell,a,q):=\sum_{\substack{\ell p\le y\\ \ell p\equiv a\!\!\!\!\!\pmod{q}}}1,
	$$
	where the symbols $p$ will always be primes. Then we have
	
	\begin{lemma}\label{lemma2}
		Let $a\neq 0$ be a given integer, and let $A>0$ and $\varepsilon>0$. For any well factorable function $\lambda(q)$ of level $Q$, the following estimate
		$$
		\sum_{\substack{q\le Q\\ (a,q)=1}}\lambda(q)\sum_{\substack{L_1\le \ell \le L_2\\ (\ell,q)=1,~2|\ell}}\bigg(\pi(x;\ell,a,q)-\frac{\pi\big(x/\ell\big)}{\varphi(q)}\bigg)\ll \frac{x}{(\log x)^A}
		$$
		holds for $Q\le x^{4/7-\varepsilon}$ and $1\le L_1\le L_2\le x^{1-\varepsilon}$, where the implied constants depend only on $a, A$ and $\varepsilon$, $\varphi(n)$ is the Euler totient function.
	\end{lemma}
	
	\begin{remark}\label{remark:well}
		In fact, the original statement of Wang's proposition  \cite[Proposition 3.2]{Wang2} is slightly different from Lemma \ref{lemma2}.  For our applications, we add the additional restriction $2|\ell$ here in Lemma \ref{lemma2}. 
		The proof is almost the same as Wang \cite[Proposition 3.2]{Wang2}.  After applying Heath-Brown's identity, we shall consider the sum
		\begin{align*}
			& \Delta(L\mid M_1, \ldots, M_j\mid N_1, \ldots, N_j;\, q,\, a)
			\\\noalign{\vskip 0,5mm}
			& \hskip 5mm
			:= \mathop{{\sum}^*}_{\substack{2\ell m_1\ldots m_jn_1\ldots n_j\equiv a ({\rm mod}\,q)\\ \ell\in \mathscr{L},\, m_i\in \mathscr{M}_i ,\, n_i\in \mathscr{N}_i}}
			\mu(m_1)\ldots \mu(m_j)
			- \frac{1}{\varphi(q)} \mathop{{\sum}^*}_{\substack{(2\ell m_1\ldots m_jn_1\ldots n_j,\, q)=1\\ \ell\in \mathscr{L},\, m_i\in \mathscr{M}_i ,\, n_i\in \mathscr{N}_i}}
			\mu(m_1)\ldots \mu(m_j),
		\end{align*}
		where $\sum^*$ means that the summation is restricted to numbers $m_1, \ldots, m_j, n_1, \ldots, n_j$ free of prime factors  $<z$ and $\mathscr{L}, \mathscr{M}_i, \mathscr{N}_i$ are intervals of the type 
		$$
		\mathscr{L} := [(1-\Delta)L,\, L[,
		\quad
		\mathscr{M}_i := [(1-\Delta)M_i,\, M_i[,
		\quad
		\mathscr{N}_i = [(1-\Delta)N_i,\, N_i[
		$$
		with
		$$
		L M_1\ldots M_jN_1\ldots N_j=x,\qquad
		\max(M_1, \ldots, M_j)< x^{1/7}
		$$
		and $\Delta=(\log x)^{-A_1}$. Here $A_1$ is a sufficiently large constant.
		In the case $L=x^{\nu_0} \geqslant x^{3/7}$, we apply Theorem 5 of \cite{BFI86} with $M=L$, and here the coefficient $2$ is attached to $m_1\cdots m_jn_1\cdots n_j$. Otherwise, we shall apply Theorems 1, 2 and $5^*$ separately according to the partial product of $ M_1, \ldots, M_j, N_1, \ldots, N_j$ is located in some given intervals, and in these cases the coefficient $2$ is attach to $\ell$.
		
		It seems that we may further generalize Lemma \ref{lemma2} to 
		\begin{align*}
			\sum_{(a,\, q)=1}\lambda(q) \sum_{\substack{L_1\leqslant \ell\leqslant L_2\\ (\ell,\, q)=1}}
			f(\ell)\bigg(\pi(x;\, \ell, a, q) -\frac{\pi(x/\ell)}{\varphi(q)}\bigg)
			\ll_{a, A, \varepsilon} \frac{x}{(\log x)^A}
		\end{align*}
		for some smooth function  $f(\ell)\ll \tau(\ell)^B$ with $B>0$. The main difference between the proofs is that we need an analogue of Theorem 5 in \cite{BFI86} with coefficient $\alpha_{\ell}\equiv 1$ replaced by smooth function $f(\ell)$, which is just Bombieri-Friedlander-Iwaniec have done in the proof of Theorem 5 in \cite{BFI86}. Here we do not pursue the details.
	\end{remark}

	The last lemma is another conjecture of Chen and Chen \cite{CC} which was later confirmed by Wu \cite[Theorem 2]{Wu}.
	\begin{lemma}\label{lemma5}
		For $0<c<1$, let 
		$$
		T'_c(x)=\#\big\{p\le x:\ P^+(p-1)\ge x^c\big\}.
		$$
		Then for sufficiently large $x$ we have
		$$
		T_c(x)=T'_c(x)+O\left(\frac{x\log\log x}{(\log x)^2}\right).
		$$
	\end{lemma}
	
	\vskip 5mm
	
	\section{Proof of Theorem \ref{thm1}}
	Now, we turn to the proof of our theorem.
	Throughout, the symbols $p$ and $p'$ will always be primes and $x$ is supposed to be sufficiently large.
	
	First, by Lemma \ref{lemma5}, it suffices to show that for any $e^{-\frac{2}{7}}< c<1$,
	$$
	\limsup_{x\rightarrow\infty}T'_c(x)/\pi(x)\le -\frac{7}{2}\log c.
	$$
	Clearly, for $c> e^{-\frac{2}{7}}>0.75$ we have
	\begin{align*}
		T'_c(x)=\sum_{\substack{p'\le x\\ P^+(p'-1)\ge x^c}}1=\sum_{x^c\le p<x}\sum_{\substack{p'\le x\\ p|p'-1}}1=\sum_{x^c\le p<x}\ \sum_{\substack{\ell p+1\le x\\ \ell p+1\,\textrm{is\, prime}\\ 2|\ell}}1
		\le \sum_{\substack{\ell\le x^{1-c}\\ 2|\ell}}\sum_{\substack{\ell p\le x\\ \ell p+1\,\textrm{is\, prime}}}1.
	\end{align*}
	We are leading to sieve out primes in the following sequence
	$$
	\mathcal{A}:=\Big\{\ell p+1: \ell\le x^{1-c}, \ell p\le x, 2|\ell\Big\},
	$$
	for some fixed $c$, $0.75<c<1$. 
	Let $\mathcal{P}=\{p:\ p\neq2\}$ and define the sieve function $S(\mathcal{A}, \mathcal{P}, z)$ by
	$$
	S(\mathcal{A}, \mathcal{P},z)=\big\{a\in \mathcal{A}: \big(a, P(z)\big)=1 \big\},
	$$
	where $z\le x^{1/2}$ is a parameter to be decided later. Then, we deduce from above notation that
	\begin{align}\label{start}
		T'_c(x)\le \sum_{\substack{a\in \mathcal{A}\\ a\,\textrm{is\, prime}}}1\le S(\mathcal{A}, \mathcal{P},z)+\pi(z)=S(\mathcal{A}, \mathcal{P},z)+O\big(x^{1/2}\big).
	\end{align}
In general, to estimate  $S(\mathcal{A}, \mathcal{P},z)$ by sieve methods, we need to compute $\mathcal{A}_d=\{a\in \mathcal{A}: d\mid a\}$ as the form
$$
\big|\mathcal{A}_d\big|=\frac{\omega(d)}{d}X+r\big(\mathcal{A}, d\big),
$$
where $X$ is a convenient close approximation to $|\mathcal{A}|$ which should be independent of $d$, $\omega(d)$ is a multiplicative function and $r\big(\mathcal{A}, d\big)$ is the error term. Unfortunately, in our question it seems  we can not extract such a  $X$ independent of $d$ from $\mathcal{A}_d$. To avoid this obstacle, we shall start from the fundamental inequality 
$$\lambda_d^- * \mathbf{1}\leqslant \mu * \mathbf{1} \leqslant \lambda_d^+ * \mathbf{1}$$
involving sieve weights $\lambda^\pm_d$, based on an idea of Fouvry and Tenenbaum \cite{FT96}. That is, we will apply Lemma \ref{lem:sieve} to give an upper bound for $S(\mathcal{A}, \mathcal{P},z)$.

First, by the M\"{o}bius inversion formula, we have
	\begin{align*}
	S(\mathcal{A}, \mathcal{P},z)=&\mathop{\sum\, \sum}_{\substack{\ell p\le x\\ \ell\le x^{1-c},\, 2|\ell}}\textbf{1}*\mu\big((\ell p+1,\, P(z))\big)\\[2mm]
	\leq &\mathop{\sum\, \sum}_{\substack{\ell p\le x\\ \ell\le x^{1-c},\, 2|\ell}}\
	\sum_{\substack{d|\ell p+1\\ d|P(z)}}\lambda^+_d
	\\[2mm]
	=&\sum_{\substack{ d|P(z)}}\lambda^+_d \mathop{\sum\, \sum}_{\substack{\ell p\le x\\ \ell\le x^{1-c},\, 2|\ell\\ \ell p+1\equiv 0\!\!\!\!\!\pmod{d}}}1.
\end{align*}
Next we approximate the inner sum over $p$ by $\pi(x/\ell)/\varphi(d)$ getting
\begin{align}\label{S(APz)}
	S(\mathcal{A}, \mathcal{P},z)
	&\leq \sum_{\substack{ d|P(z)}}\lambda^+_d \sum_{\substack{\ell\le x^{1-c}\\  2|\ell,\, (\ell,d)=1}}\frac{\pi(x/\ell)}{\varphi(d)}+\sum_{\substack{ d|P(z)}}\lambda^+_d \,r\big(\mathcal{A}, d\big)
	\nonumber\\[2mm]
	&=: M(\mathcal{A}, \mathcal{P},z)+\sum_{\substack{ d|P(z)}}\lambda^+_d \,r\big(\mathcal{A}, d\big),
\end{align}
	say, where $r(\mathcal{A}, d)$ is the error term defined by
\begin{align*}
	r(\mathcal{A}, d)=\sum_{\substack{\ell\le x^{1-c}\\  2|\ell,\, (\ell,d)=1}}\bigg(\pi(x;\ell,-1,d)-\frac{\pi(x/\ell)}{\varphi(d)}\bigg).
\end{align*}

	Thanks to the arguments of Iwaniec \cite[Theorems 1 and 4]{Iwaniec2}, the sum $\sum_{ d}\lambda^+_d \,r(\mathcal{A}, d)$ can be rearranged and transformed to the following flexible form:
	$$
	\sum_{h<\exp(8/\varepsilon^2)}\sum_{d|P(z)}\lambda^+(d,h)r\big(\mathcal{A}, d\big),
	$$
	where the coefficients $\lambda^+(d,h)$ satisfy $\big|\lambda_h^{+}(d,h)\big| \le 1$ and vanish for $d>D$. Especially, $\lambda^+(d,h)$ are  well factorable of  level $D$. Then by applying  Lemma \ref{lemma2} with $D=x^{4/7-\varepsilon}$ and $L_1=1$, $L_2=x^{1-c}$,  we obtain for the error term 
\begin{align}\label{error}
	\sum_{\substack{ d|P(z)}}\lambda^+_d \,r(\mathcal{A}, d) &=\sum_{h<\exp(8/\varepsilon^2)}\sum_{\substack{d\le x^{4/7-\varepsilon}}}\lambda^+(d,h)r\big(\mathcal{A}, d\big)\nonumber\\
	&=\sum_{h<\exp(8/\varepsilon^2)}\sum_{\substack{d\le x^{4/7-\varepsilon}}}\lambda^+(d,h)\sum_{\substack{\ell\le x^{1-c}\\  2|\ell, (\ell,d)=1}}\bigg(\pi(x;\ell,-1,d)-\frac{\pi(x/\ell)}{\varphi(d)}\bigg)\nonumber\\
	&\ll_{\varepsilon} \frac{x}{(\log x)^A}
\end{align}
for any $A>0$, which is negligible.

Now it remains to estimate the main term $M(\mathcal{A}, \mathcal{P},z)$. Changing
the order of summation, we find that
\begin{align*}
M(\mathcal{A}, \mathcal{P},z)=&
\sum_{\substack{\ell\le x^{1-c}\\  2|\ell}}\pi(x/\ell) \sum_{\substack{ d|P(z)\\ (d,\,\ell)=1}}\frac{\lambda^+_d}{\varphi(d)}
\\[2mm]
=&
\sum_{\substack{\ell\le x^{1-c}\\  2|\ell}}\pi(x/\ell) \sum_{\substack{ d|P'(z)}}\frac{\lambda^+_d}{\varphi(d)}
\end{align*}
where $$P'(z)=\prod_{p\leq z,\, p\in \mathscr{P}}p, \quad \mathscr{P}=\big\{p:\ p\nmid \ell\big\}.$$
Then we apply Lemma \ref{lem:sieve} with $\omega(p)=p/\varphi(p),\, z=D=x^{4/7-\varepsilon}$ and some large constant $L$ getting
\begin{align*}
	M(\mathcal{A}, \mathcal{P},z)& \leq \Big\{F(1)+o(1)\Big\} \sum_{\substack{\ell \leq x^{1-c}\\ 2|\ell }}\frac{x}{\ell\log(x/\ell)}
	\ \prod_{\substack{2<p\leq z\\ p\nmid \ell }}
	\Big(1-\frac{1}{p-1}\Big)
	\\\noalign{\vskip 1mm}
	& =\Big\{F(1)+o(1)\Big\}\,\prod_{\substack{2<p\leq z}}
	\Big(\frac{p-2}{p-1}\Big) \sum_{\substack{\ell \leq x^{1-c}\\ 2|\ell }}\frac{x}{\ell\log(x/\ell)}\prod_{2<p\leq z,\, p|\ell}\Big(\frac{p-1}{p-2}\Big)
	\\\noalign{\vskip 1mm}
	& =\Big\{F(1)+o(1)\Big\}\,\prod_{\substack{2<p\leq z}}
	\Big(\frac{p-2}{p-1}\Big) \sum_{\substack{\ell \leq x^{1-c} }}\frac{x}{2\ell\log(x/2\ell)}\prod_{2<p\leq z,\, p|2\ell}\Big(\frac{p-1}{p-2}\Big)
\end{align*}
Considering 
$$
\frac{1}{\log(x/2\ell)}=\frac{1}{\log(x/\ell)}\Big\{1+O\Big(\frac{1}{\log x}\Big)\Big\}
\quad {\rm and}\quad 
\prod_{2<p\leq z,\, p|\ell}\Big(\frac{p-1}{p-2}\Big)=\prod_{p>2,\, p|\ell}\Big(\frac{p-1}{p-2}\Big)
$$
since $z=x^{4/7-\varepsilon}>x^{1-c}\geq\ell$, we then arrive at
\begin{align}\label{M(APz)}
	M(\mathcal{A}, \mathcal{P},z)& =x\Big\{\frac{F(1)}{2}+o(1)\Big\}\,\prod_{\substack{2<p\leq z}}
	\Big(\frac{p-2}{p-1}\Big) \sum_{\substack{\ell \leq x^{1-c}}}\frac{1}{\ell\log(x/\ell)}\prod_{p>2,\, p|\ell}\Big(\frac{p-1}{p-2}\Big).
\end{align}
Next we define the multiplicative function $H(\ell)$  by
$$
H(\ell):=\prod_{p>2,\, p|\ell}\Big(\frac{p-1}{p-2}\Big).
$$
The partial sum of $H(\ell)$ is well studied by La Bret\`{e}che, Pomerance and Tenenbaum \cite{BPT05} by employing the Selberg-Delange method: for any nonnegative integer $v$,
$$
\sum_{\substack{\ell\leq y\\ (\ell,\, a)=1}}H(2^v\ell) \sim \frac{Ay}{2^{\varepsilon(a)}H(a)}\qquad (y\to\infty),
$$
where $\varepsilon(a)=1$ if $2\mid a$ and $\varepsilon(a)=0$ if $2\nmid a$, and $A$ is an absolute constant defined by
$$
A:= \prod_{p>2}\Big(1+\frac{1}{p(p-2)}\Big).
$$
By taking $a=1$ in the above asymptotic formula and by partial summation, we infer that
\begin{align}\label{sum-ell}
\sum_{\substack{\ell \leq x^{1-c}}}\frac{1}{\ell\log(x/\ell)}H(\ell)=
 \int^{x^{1-c}}_{1}\frac{{A\rm{d}}t}{t\log(x/t)}+o(1)
 =  A\log\frac{1}{c} +o(1). 
\end{align}
Inserting \eqref{sum-ell} into \eqref{M(APz)} we get
\begin{align*}
	M(\mathcal{A}, \mathcal{P},z)& \leq x\Big\{\frac{F(1)}{2}+o(1)\Big\}A\log\frac{1}{c}\,
	\prod_{\substack{2<p\leq z}}\Big(\frac{p-2}{p-1}\Big).
\end{align*}
Note that $F(1)=2\textrm{e}^{\gamma}$, and by Merten's theorem
\begin{align*}
	\prod_{\substack{2<p\leq z}}\Big(\frac{p-2}{p-1}\Big)
	&=\prod_{\substack{2<p\leq x^{4/7-\varepsilon}}}
	\Big(\frac{p-2}{p-1}\Big)\Big(1-\frac{1}{p}\Big)^{-1}\prod_{\substack{2<p\leq x^{4/7-\varepsilon}}}\Big(1-\frac{1}{p}\Big)
	\\\noalign{\vskip 0,5mm}
	&=\frac{2\textrm{e}^{-\gamma}+o(1)}{\log x^{4/7-\varepsilon}}\prod_{\substack{2<p\leq x^{4/7-\varepsilon} }}\frac{p(p-2)}{(p-1)^2}
	\\\noalign{\vskip 0,5mm}
	&=\frac{\frac{7}{2}\textrm{e}^{-\gamma}+o(1)}{\log x}\prod_{\substack{p>2 }}\Big(1+\frac{1}{p(p-2)}\Big)^{-1}\Big\{1+O\Big(\frac{1}{x}\Big)
	\Big\}
	\\\noalign{\vskip 0,5mm}
	&=\frac{1}{\log x}
	\Big(\frac{7}{2}\textrm{e}^{-\gamma}A^{-1}+o(1)\Big),
\end{align*}
therefore we find
\begin{align}\label{M(APz)}
	M(\mathcal{A}, \mathcal{P},z)& \leq \Big(-\frac{7}{2}\log c+o(1)\Big)\frac{x}{\log x}.
\end{align}

Finally, combining \eqref{S(APz)}, \eqref{error} and \eqref{M(APz)}, we conclude that
	\begin{align*}
		S(\mathcal{A}, \mathcal{P},z)\le \Big(-\frac{7}{2}\log c+o(1)\Big)\frac{x}{\log x},
	\end{align*}
and hence, $T'_c(x)$ is estimated from (\ref{start}) as follows
	\begin{align*}
		T'_c(x)\le S(\mathcal{A}, \mathcal{P},z) +O(x^{1/2})\le
		\Big(-\frac{7}{2}\log c+o(1)\Big)\frac{x}{\log x},
	\end{align*}
	whence
	\begin{align*}
		\limsup_{x\rightarrow\infty }T'_c(x)/\pi(x)\le -\frac{7}{2}\log c
	\end{align*}
	for any $e^{-\frac{2}{7}}< c<1$. This completes the proof of Theorem \ref{thm1}.

	\section*{Acknowledgments} 
	The authors are also very grateful to Zhiyuan Yang
	for pointing out a technical mistake in the first version of the
	paper.
	
	This work was supported by the National Key Research and Development Program of China (Grant No. 2021YFA1000700).
	
	Y. C. was supported by National Natural Science Foundation of China  (Grant No. 12201544) and China Postdoctoral Science Foundation (Grant No. 2022M710121).
	
	Z. W. was supported by National Natural Science Foundation of China  (Grant No. 12322101) and National Natural Science Foundation of Shandong Province  (Grant No. ZR2023YQ004).

\end{document}